\documentclass[11pt]{llncs}
\usepackage{amssymb,amsmath}
\usepackage{epsfig}
\usepackage{epstopdf}
\usepackage{graphicx}
\usepackage{llncs.cls}
\newtheorem{theorems}{Theorem}
\newtheorem{defn}{Definition}

\usepackage{graphicx}
\usepackage{caption}
\title{A New Characterisation of Total Graphs}
\author{Ravi Goyal,  Mahipal Jadeja, Rahul Muthu}
\institute{Dhirubhai Ambani Institute of Information  Communication Technology, Gandhinagar, India\\ \email{\{ravi\_goyal, jadeja\_mahipal, rahul\_muthu\}@daiict.ac.in}}

\date{}
\begin{document}
\maketitle

\begin{abstract} Graphs constructed to translate some graph problem into another graph problem are usually called auxiliary graphs. Specifically, total graphs of simple graphs are used to translate the total colouring problem of the original graph into a vertex colouring problem of the transformed graph.

In this paper, we obtain a new characterisation of total graphs of simple graphs. We also design algorithms to compute the inverse total graphs when the input graph is a total graph. These results improve over the work of Behzad, by using novel observations on the properties of the local structure in the neighbourhood of each vertex. The earlier algorithm was based on {\sc bfs} and distances. Our theorems result in partitioning the vertex set of the total graph into the original graph and the line graph efficiently.

We obtain constructive results for special classes, most notably for total graph of complete graphs.

\end{abstract}

{\bf keywords}: Total colouring$\cdot$Edge colouring$\cdot$Line graphs$\cdot$\\Total graphs$\cdot$Inverse total graphs
\vspace{-0.3cm}
\section{Introduction}\label{SecIntro}
\vspace{-0.3cm}
A total colouring of a simple graph is a simultaneous assignment of labels to its vertices and edges such that adjacent vertices get distinct colours, adjacent edges get distinct colours and the colour of each edge is distinct from the colours of its endpoint vertices (or equivalently the colour of each vertex is distinct from the colours of its incident edges). It is thus a combination of a proper vertex colouring, a proper edge colouring and a further restriction on the interplay between these colourings. The notion of total colouring was introduced by Behzad [2] and Vizing [8] and those papers also conjectured that $\chi_T(G)\le\Delta(G)+2$. It is immediate that $\chi_T(G)\ge\Delta(G)+1$, since a vertex of maximum degree and its incident edges must all get distinct colours. A lot of work has been done on total colouring, based on frugal colouring, the list colouring conjecture etc.

In this work our aim is not to make headway in the total colouring problem but rather to get a complete characterisation of the class of total graphs of simple graphs. Towards this end, we of course need a precise definition of what we mean by the total graph of a given graph, which is an existing concept. This also brings along the concept of line graphs. The vertex set of the total graph of a simple graph can be partitioned into two sets, one corresponding to the vertex set of the original graph (inverse total graph) and the other the line graph of the original graph, with crossing edges between these two vertex sets. A result characterising total graphs was obtained by Behzad [1] in 1970, and earlier works obtained interesting properties on this class of graphs [4]. The fact that any total graph has a unique preimage under the inverse total graph operator was proved by [3].

We obtain a new characterisation of total graphs based on the induced subgraphs on the neighbourhood of maximum degree vertices. These characterisations allow us to distinguish vertex vertices from edge vertices among the vertices of maximum degree. Using this characterisation, we develop an efficient algorithm which iteratively creates the partition of the vertex set of the candidate total graph into its inverse total graph and line graph. 

The paper is organised as follows. Definitions and notation used in the paper are presented in Section~\ref{SecDefnNot}. Basic properties and theorems on total graphs are presented in Section~\ref{SecTotGraIntro}. Section~\ref{Sectotalcomplete} presents our results for special classes, the most elegant among which is the explicit construction for complete graphs. In Section~\ref{SecTotGraChar} we present our theorems characterising the two classes of vertices and develop our theorem into an algorithm for reconstructing the inverse total graph of a given total graph. We summarise our work and indicate possible future directions for research in Section~\ref{SecConcl}. 
\vspace{-0.3cm}
\section{Definitions and Notation}\label{SecDefnNot}
\vspace{-0.3cm}
In this section we present the basic definitions and notation we use throughout the paper. All graphs we consider are finite, simple, undirected and connected. This problem for graphs with more than one component are strsaightforward extensions of the connected case.
\vspace{-0.1cm} 
\begin{defn}\label{deflngr}
The {\bf line graph} $L(G)$ of a graph $G=(V,E)$ is defined as the graph with vertex set having one vertex corresponding to each edge in $G$ and an edge between two vertices of $L(G)$ precisely when the edges of $G$ that those vertices correspond to, have a common endpoint. 
\end{defn}
Not all graphs are line graphs of a simple graph (or for that matter even of a multigraph). The line graph of a graph is however well defined, and results in a unique graph. The concept of line graph can thus be viewed as a non-surjective function from the set of simple graphs to the set of simple graphs. When we limit attention to connected simple graphs, it is further known that this function is almost injective. There is only one pair of distinct graphs with the same line graph namely the three vertex cycle $C_3$ and the three leaf star $K_{1,3}$ both of whose line graphs are $C_3$. 

Thus it is an interesting question both combinatorially and algorithmically as to which graphs are line graphs of some simple graph, and computing the inverse-line graph (the preimage under the line graph function) for those graphs which are line graphs [6] Iterating the line graph operator and its inverse, as well as properties of the resulting sequence of graphs has also been studied extensively in the literature  [5][7]

An obvious application of the concept of line graphs is to translate edge colouring problems into equivalent vertex colouring problems on the line graph. Similarly, we define the total graph of a graph. Again, this transformation renders a total colouring problem on a graph into a vertex colouring problem on an equivalent graph. The total graph concept can also be viewed as a function from the class of graphs to the class of graphs. It has been established that the total graph viewed as a function from the class of graphs to the class of graphs is injective. 

\begin{defn}\label{deftotgr}
The total graph $T(G)$ of a graph $G=(V,E)$ has as vertex set one vertex for each edge and vertex in $G$. Two vertices in $T(G)$ are adjacent precisely when the elements (vertex or edge) of $G$ they represent are adjacent/incident to each other in $G$.
\end{defn}

\begin{defn}\label{defnvertcatg}
Thus the vertex set of the total graph of a graph can be partitioned into:
\begin{enumerate}
\item  The vertices of the original graph $($we call such a vertex a {\bf vertex vertex}$)$ 
\item The vertices of the line graph $($we call such a vertex an {\bf edge vertex}$)$.
\end{enumerate}
\end{defn}
\begin{defn}\label{defmixclq}
A {\bf mixed clique} in a total graph is a clique which has at least one vertex from the set of vertex vertices and at least one vertex from the set of edge vertices in a valid partition of the total graph into the inverse total graph and its line graph.
\end{defn}
\begin{defn}\label{defnpureclq}A {\bf pure clique} in a total graph is a clique consisting exclusively of vertex vertices or exclusively of edge vertices.
\end{defn}
\vspace{-0.3cm}
\section{Total Graphs: Basic Properties}\label{SecTotGraIntro}
\vspace{-0.3cm}
The vertex set of a total graph $H=T(G)$ has a partition into two disjoint sets, such that the induced subgraph on one part is the inverse total graph $G$ and the other is the line graph $L(G)$. In several total graphs, this partitioning is non-unique, but the different partitions are all isomporphic to each other. These are called the vertex part and the edge part respectively. The bipartite subgraph induced by the division as above is 2-regular on vertices in the edge part. This is because each edge in the inverse total graph has precisely two endpoints. Thus an algorithm for computing the inverse total graph of a given graph (if it is a total graph) is based on finding such a division or partition of the vertex set.
\begin{theorems}\label{ThmMixedClq}
The largest mixed clique consisting of at least two vertex vertices is in a total graph is of size 3.  
\end{theorems}
\vspace{-0.3cm}
\begin{proof}
Consider a clique consisting of three vertex vertices. Clearly there is no edge in any graph incident to three distinct vertices. Hence this clique cannot be augmented to include any edge vertices, and thus is not the subset of any mixed clique.
It follows that a mixed clique can consist of at most two adjacent vertex vertices. In this case this can be augmented by only one vertex, the edge vertex representing the link between these two adjacent vertex vertices. 
 \vspace{-.7cm}\begin{flushright}$\Box$\end{flushright}
\end{proof}
\vspace{-0.3cm}
Thus a maximal mixed clique is either:
\vspace{-0.3cm}
\begin{itemize}
\item A vertex vertex and all its adjacent edge vertices; or
\item Two adjacent vertex vertices and the connecting edge vertex.
\end{itemize}
In summary, a maximal mixed clique is either of size 3 or of size $k+1$ where $2k$ is the degree of its only vertex vertex in the total graph.
\vspace{-0.3cm}
\subsection{Vertex Degrees of $H=T(G)$ in terms of Vertex Degrees of $G$.}\label{SubSecdegseq}
\vspace{-0.3cm}
In this subsection, we derive results expressing the degree sequence of a total graph in terms of the degree sequence of its inverse total graph.
\begin{theorems}\label{ThmDegVertVert}
The degree of a vertex vertex in a total graph is 2 times the degree of the original vertex in the inverse total graph.
\end{theorems}
\vspace{-0.3cm}
\begin{proof}
In the total graph a vertex vertex is adjacent to vertices corresponding to its original neighbours as well as its incident edges in the original graph. See Figure~\ref{fig:vertexedge}. 
 \vspace{-.7cm}\begin{flushright}$\Box$\end{flushright}
\end{proof}
\vspace{-0.3cm}
\begin{theorems}\label{ThmDegEdgeVert}
The degree of an edge vertex in a total graph is equal to the sum of the degrees of the endpoint vertices of the original edge in the inverse total graph.
\end{theorems}
\vspace{-0.3cm}
\begin{proof}
In the total graph an edge vertex is adjacent to its two endpoints (which are both vertex vertices), and the other edges incident to these endpoints (which are all edge vertices). Thus if its endpoints are $u$ and $v$, its degree is $2+(d_u-1)+(d_v-1)=d_u+d_v$. See Figure~\ref{fig:vertexedge}. 
 \vspace{-.7cm}\begin{flushright}$\Box$\end{flushright}
\end{proof}
\vspace{-0.3cm}
\begin{theorems}\label{ThmDegTotReg}
The total graph of a graph is regular if and only if the original graph is regular.
\end{theorems}
\vspace{-0.3cm}
\begin{proof}
This is an immediate consequence of Theorems~\ref{ThmDegVertVert} and \ref{ThmDegEdgeVert}.
 \vspace{-.7cm}\begin{flushright}$\Box$\end{flushright}
\end{proof}
\vspace{-0.3cm}
\begin{theorems}\label{cliqueververedge}
Consider a maximal mixed clique of size 3 in a total graph. Let the two vertex vertices be of degrees $a$ and $b$, with $a>b$. Then the degree of the edge vertex of this clique is $c=\frac{a+b}{2}$. Clearly $a>c>b$.
\end{theorems}
\vspace{-0.3cm}
\begin{proof}
This follows from Definition~\ref{defmixclq} and Theorem~\ref{ThmMixedClq}
 \vspace{-.7cm}\begin{flushright}$\Box$\end{flushright}
\end{proof}
\vspace{-0.3cm}
\begin{theorems}
Given a vertex vertex $v_a$ of degree $2k$ in a total graph, its neighbours can be divided into $k$ pairs representing:
\begin{itemize}
\item  Its distinct incident edges $(v_1,\ldots,v_k)$ in the inverse total graph and
\item the other endpoints of those edges $v_1',\ldots,v_k')$ respectively. 
\end{itemize}
The $k$ pairs are $\{(v_1,v'_1),\ldots,(v_k,v_k')\}$. 
\end{theorems}
\vspace{-0.3cm}
\begin{proof}
See Figure~\ref{fig:vertexedge}. The degrees of each pair is related to the degree of the selected vertex vertex according to Theorem~\ref{cliqueververedge}
 \vspace{-.7cm}\begin{flushright}$\Box$\end{flushright}
\end{proof}
\vspace{-0.3cm}
\subsection{Relationship between the Number of Vertices and Edges of the Original Graph and its Total Graph
}\label{SubSecRelation}
\vspace{-0.3cm}
Here we extend the work of the previous subsection in a natural way obtaining results expressing the number of vertices and edges of a total graph in terms of the corresponding parameters of the inverse total graph.
\begin{theorems}\label{relation}
Let $T(G)=(V',E')$ for the given $G=(V,E)$.

\begin{enumerate}
	\item $|V'|=|V|+|E|$
\item $|E'| \le  |E|(|V|+1)$
\end{enumerate}
\end{theorems}
\vspace{-0.3cm}
\begin{proof}
\begin{itemize}
\item The total graph of a given graph contains the original graph and the line graph of the original graph as disjoint induced subgraphs spanning all its vertices. The line graph of the original graph G has $|E|$ vertices. Therefore, $|V'|=|V|+|E|$

\item The edge set of the total graph can be partitioned into $3$ disjoint sets:

\begin{enumerate}
\item[]  {\bf A:} Edge set of the original graph containg $|E|$ edges.
\item[] {\bf B:} Edges which are present between the original graph and its line graph. There are  $2|E|$ such edges.

Note: The line graph of the given $G$ has $|E|$ vertices and each vertex of the line graph is connected to exactly $2$ vertices of the original graph. 
 \item[] {\bf C}: Edges which are present within the line graph. The number of such edges is \[\left(\Sigma_{{v_i}\in V(L(G))}\frac{(d(v_i))^2}{2}\right)- |E|\]
\end{enumerate}
\item Explanation for edges present within the line graph:

The number of edges in the subgraph induced by any subset $X$ of vertices of a graph $G$ is given by \[(\Sigma_{x\in X} d_G(x))-|E(X,\overline{X})|\] 
Applying this rule, gives us the claimed result for portion {\bf C}.\\

Now the maximum possible degree for any vertex of the given simple connected graph is $|V|-1$. \\
 Thus, 

$\begin{array}{rcl}d(v_1)d(v_1)+\cdots+d(v_n)d(v_n)&\le &(|V|-1) \times (\mbox{degree sum})\\
 &\le &(|V|-1)  2|E|
\end{array}$

\item  Now the total number of edges in $T(G)= |A \cup B \cup C| $\\

$\begin{array}{rcl}
\therefore |E'| & \le  &|E|+2|E|+  (|V|-1) |E|- |E|\\
& \le & 2|E|+  (|V|-1) |E|\\
& \le & |E|+  |V| |E|\\
& \le & |E|(|V|+1)\\
\end{array}$
\end{itemize}
 \vspace{-.7cm}\begin{flushright}$\Box$\end{flushright}
\end{proof}
\vspace{-0.3cm}
For instance $T(K_4\setminus {e})$ has 9 vertices and 23 edges. These numbers conform to the constraints derived by us above.
\vspace{-0.3cm}
\section{Results on Paths, Cycles and Complete Graphs}\label{Sectotalcomplete}
\begin{figure}[h]
\centering
  \begin{tabular}{@{}cccc@{}}
    \includegraphics [height= 4cm, width=.50\textwidth]{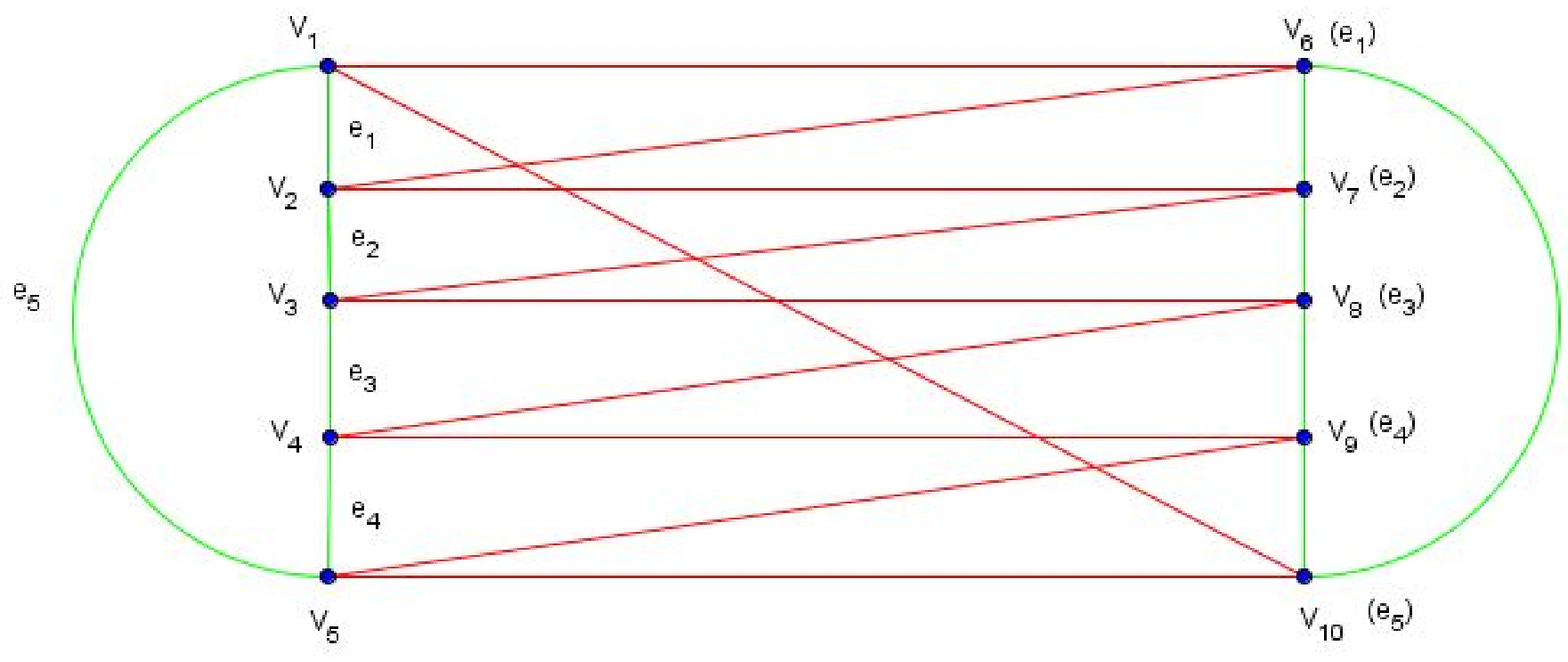} &
    \multicolumn{2}{c}{\includegraphics[height= 4cm, width=.50\textwidth]{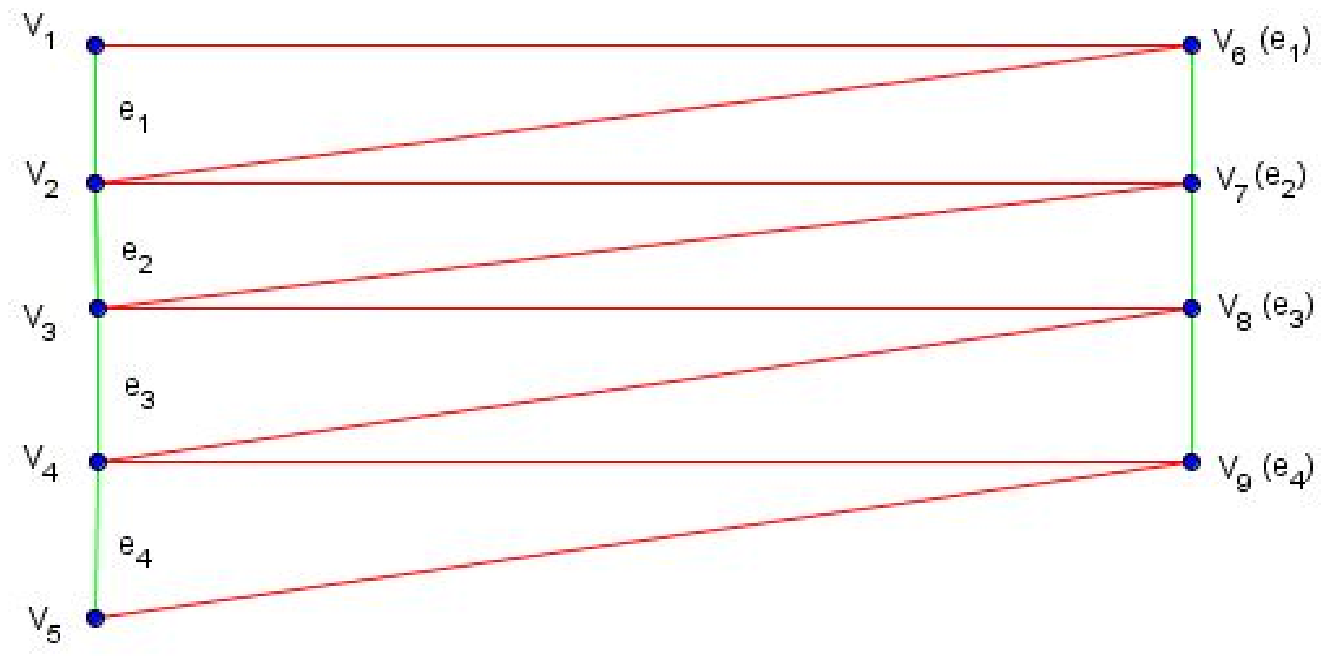}}
  \end{tabular}
  \caption{Example of hamiltonian cycle and hamiltonian path for $T(C_5)$ and $T(P_5)$ \label{hamiltonian}}
\end{figure}

\vspace{-0.3cm}
We give explicit constructions for the total graphs of paths, cycles and complete graphs in this section. The first result is for cycles.
\begin{theorems}\label{ThmCycles}
A graph is the total graph of the cycle $C_n$ on $n$ vertices if and only if:
\begin{itemize}
\item  it has $2n$ vertices $\{v_1,\ldots,v_{2n}\}$
\item  Contains two vertex disjoint cycles of length $n$ each, labeled \[\mathcal{C}_1=\{v_1,\ldots,v_n\}\mbox{ and }\mathcal{C}_2=\{v_{n+1},\ldots,v_{2n}\}\]
\item  Contains a hamiltonian cycle edge disjoint with the afore mentioned cycles with the order of vertices being \[\mathcal{C}_{ham}=\{v_1,v_{n+1},v_2,v_{n+2},\ldots,v_n,v_{2n},v_1\}\].
\end{itemize}
\end{theorems}
\vspace{-0.3cm}
\begin{proof} This result is self explanatory. See ~Figure \ref{hamiltonian}.
\vspace{-.7cm}\begin{flushright}$\Box$\end{flushright}
\end{proof}
\vspace{-0.3cm}

Now we give the structure of the total graph of paths.
\begin{theorems}\label{ThmPaths}
A graph is the total graph of the path $P_n$ if and only if:
\begin{itemize}
\item  it has $2n-1$ vertices $\{v_1,\ldots,v_{2n-1}\}$
\item disjoint paths \[\{v_1,\ldots,v_n\}\mbox{ of length }n-1\mbox{ and} \]\[\{v_{n+1},\ldots,v_{2n-1}\}\mbox{ of length }n-2\]
\item  a hamiltonian path edge disjoint with the afore-mentioned paths with order of vertices being $\{v_1,v_{n+1},v_2,v_{n+2},\ldots,v_{n-1},v_{2n-1},v_n\}$.
\end{itemize}
\end{theorems}
\vspace{-0.3cm}
\begin{proof}
This result is self explanatory. See ~Figure \ref{hamiltonian}.
\vspace{-.72cm}\begin{flushright}$\Box$\end{flushright}
\end{proof}
\vspace{-0.3cm}
The following theorems give the structure and other parameters of the total graphs of complete graphs.
\begin{theorems}\label{ThmCompNumb}
Let $T(G)=(V',E')$ for the given $K_{|V|}$.
\begin{enumerate}
	\item $|V'|= \frac{|V||(|V|+1)}  {2}$
\item $|E'| = \frac {|V|(|V|-1)(|V|+1)} {2}$
\item $\forall v\in K_{|V|}, d_v=2(|V|-1)$.
\end{enumerate}
\end{theorems}
\vspace{-0.3cm}
\begin{proof} This follows from Theorems \ref{relation} and \ref{ThmDegVertVert}. 
\vspace{-.7cm}
\end{proof}
\vspace{-0.3cm}
\begin{theorems}\label{ThmCompLineTot}
$L(K_{n})=T(K_{n-1})$ where $L(K_{n})$ denotes the line graph of $K_{n}$.
\end{theorems}
\vspace{-0.3cm}
\begin{proof}
\begin{itemize}
\item The line graph of $K_{n}$ has exactly $n$ cliques each of size $n-1$. These are formed by the  stars induced by the edges incident to each of the $n$ vertices of $K_n$. Each pair of these cliques share a unique common vertex. For instance the first clique has a distinct vertex common with each of the remaining $n-1$ cliques. The same holds for each of these cliques.
\item $T(K_{n-1})$ contains $L(K_{n-1})$.  It follows that there are exactly $n-1$ maximal cliques each of size $n-2$ in its induced line graph. 
\item $T(K_{n-1})$ also contains a copy of $K_{n-1}$, ignoring the vertices of the line graph portion. This is basically the inverse total graph. 
\item Thus, in total $n$ cliques are present, $n-1$ of which are of size $n-2$ and the remaining $1$ of size $n-1$. This follows from the previous two points.
\item Each clique of size $n-2$ present in line graph is adjacent to exactly $1$ vertex of the original graph because the clique is formed by the vertices corresponding to  edges incident on the single vertex of the original graph. 
\item Hence it is possible to include exactly $1$ more vertex in each of the $n-1$ cliques of size  $n-2$ present in the line graph.
\item We conclude that a total of $n$ cliques each of size $n-1$ are present in $T(K_{n-1})$. This is the same as in $L(K_{n})$ including the overlapping pattern among these cliques.
\end{itemize}
\vspace{-.7cm}\begin{flushright}$\Box$\end{flushright}
\end{proof}
\vspace{-0.1cm}
{\bf \underline{Algorithm 1: Is a given graph the total graph of $K_n$}}
\begin{enumerate}
\item Check if the input graph satisfies the conditions of Theorem~\ref{ThmCompNumb} and compute $n$ if yes. Else reject the input. 
\item If $n=2\mbox{ or }3$ check directly for isomorphism with $L(K_3)$ or $L(K_4)$. 
\item Else look for a clique of size 4. From Theorem~\ref{ThmMixedClq} this is a pure clique. Extend it greedily to a maximal clique of size $n$. This must be possible if the given graph is $L(K_n)$. 
\item Check if the subgraph induced by the remaining vertices is \\$T(K_{n-1})=L(K_n)$ (by Theorem~\ref{ThmCompLineTot}). Check that the cross connections are consistent between the vertex part and the line graph part. Report whether the graph is $T(K_n)$ accordingly.
\end{enumerate}
\vspace{-0.3cm}
\subsection{Direct Method for Construction of Total Graph of a Complete Graph ($T(K_n)$)} 
\begin{figure}[h]
\centering
  \begin{tabular}{@{}cccc@{}}
    \includegraphics[width=.60\textwidth]{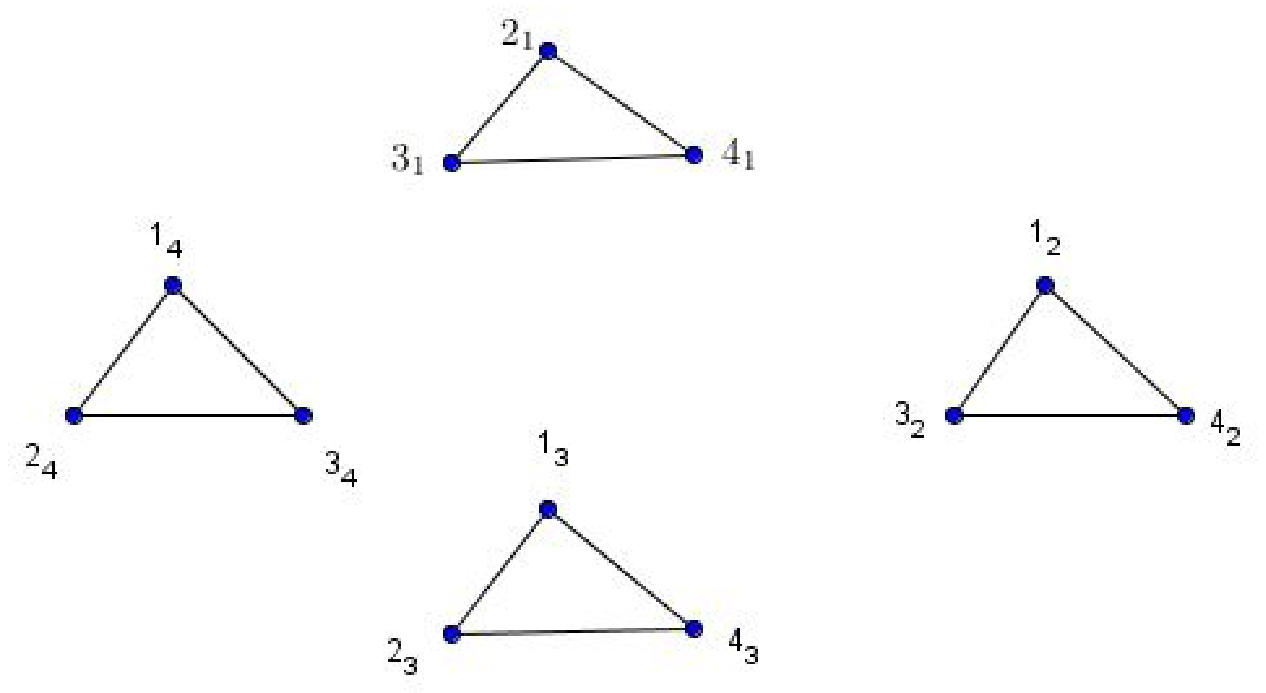} &
    \multicolumn{2}{c}{\includegraphics[width=.50\textwidth]{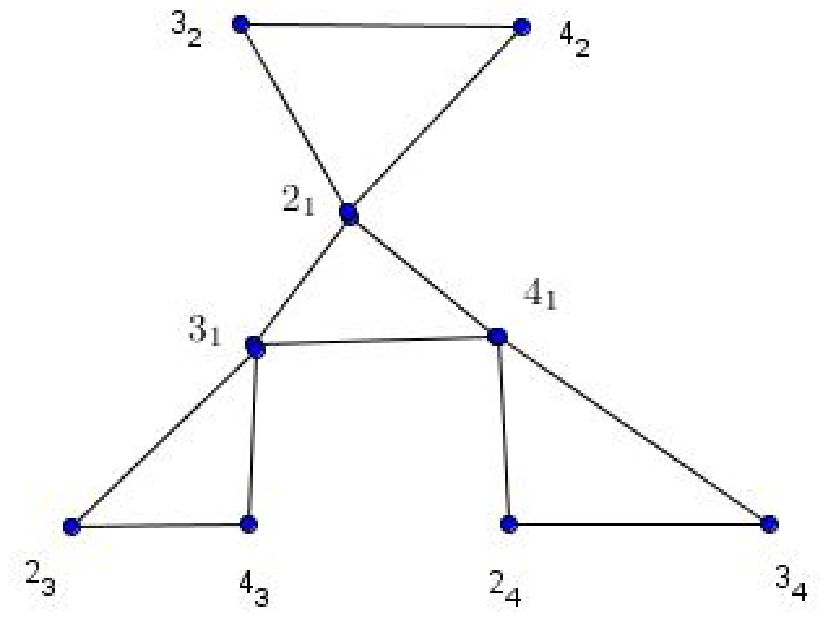}}
  \end{tabular}
  \caption{Initial steps of direct construction of the total graph of $K_{3}$\label{fig:construction}.  Group numbers are subscripts;  vertex number within each group is in normal font.}
	
\end{figure}

\begin{figure}[h]
\centering
  \begin{tabular}{@{}cccc@{}}
    \includegraphics[height=3.5cm,width=.990\textwidth]{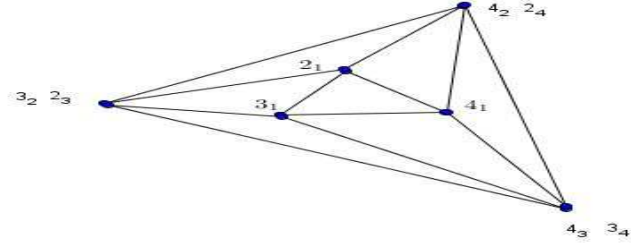} &
    \multicolumn{2}{c}{}
  \end{tabular}
  \caption{Resultant Total Graph of $K_{3}$ after step $4$ \label{fig:Resultant}}
\end{figure}

\vspace{-0.3cm}
\begin{enumerate}
\item Consider $n+1$ distinct groups, each containing exactly $n$ elements. (These correspond to the $n+1$ cliques of size $n$  of the total graph of the $K_n$).
\item Elements present in  $i^{th}$ group $G_i = \{1,\cdots, n+1\}\setminus \{i\}$. Therefore,  $|G_i|=n$ 
\item For each $G_i$, construct $K_n$ by connecting all the $n$ elements present in the group to each other. 
\item Combine the  $j^{th}$ vertex of group $i$ and the $i^{th}$ vertex of group $j$. into a single vertex. The neighbourhood of the new vertex is the union of the individual neighbourhoods of the two combined vertices. The degree of each new vertex is exactly $2(n-1)$ because degree of each original vertex is $n-1$
\item The resultant graph is the total graph of the complete graph $K_n$.
\end{enumerate}
\vspace{-0.3cm}
{\bf Notes:}\begin{itemize}
\item  No edge is destroyed during the entire procedure. 
\item Each clique has $n \choose 2$ edges and initially $n+1$ distinct cliques are considered. 
\item Thus, $n \choose 2$$(n+1)$ edges are present in the graph which remains constant through the course of this construction.
\end{itemize}
\vspace{-0.2cm}
For $T(K_3)$ this  is shown in Figures~\ref{fig:construction} and ~\ref{fig:Resultant}.  


\vspace{-0.3cm}
\section{Characterisation of Total Graphs and computing the inverse total graph}\label{SecTotGraChar}
\vspace{-0.3cm}
In this section, we prove two theorems which give conditions for a maximum degree vertex in a candidate total graph to be a vertex vertex or an edge vertex. These theorems are used to iteratively find a maximum degree vertex vertex (one is guaranteed to exist by theorem~\ref{cliqueververedge}) and partition its neighbours into vertex vertices and edge vertices. At each round a maximum degree vertex vertex is selected as a part of the inverse total graph, and is eliminated along with its edge vertex neighbours to get a smaller graph to recurse on. A partition of the vertex set of the given graph into the inverse total graph and the line graph of the inverse total graph is created iteratively, if one exists; or we infer that no such partition exists if there is a violation of the combinatorial conditions at any iteration. These theorems work only for total graphs of non-complete graphs. Thus the algorithm developed in this section also uses the algorithm of Section~\ref{Sectotalcomplete}, when appropriate, to handle the case of complete graphs.

\begin{figure}[h]
\centering
  \begin{tabular}{@{}cccc@{}}
    \includegraphics [height= 5cm, width=.50\textwidth]{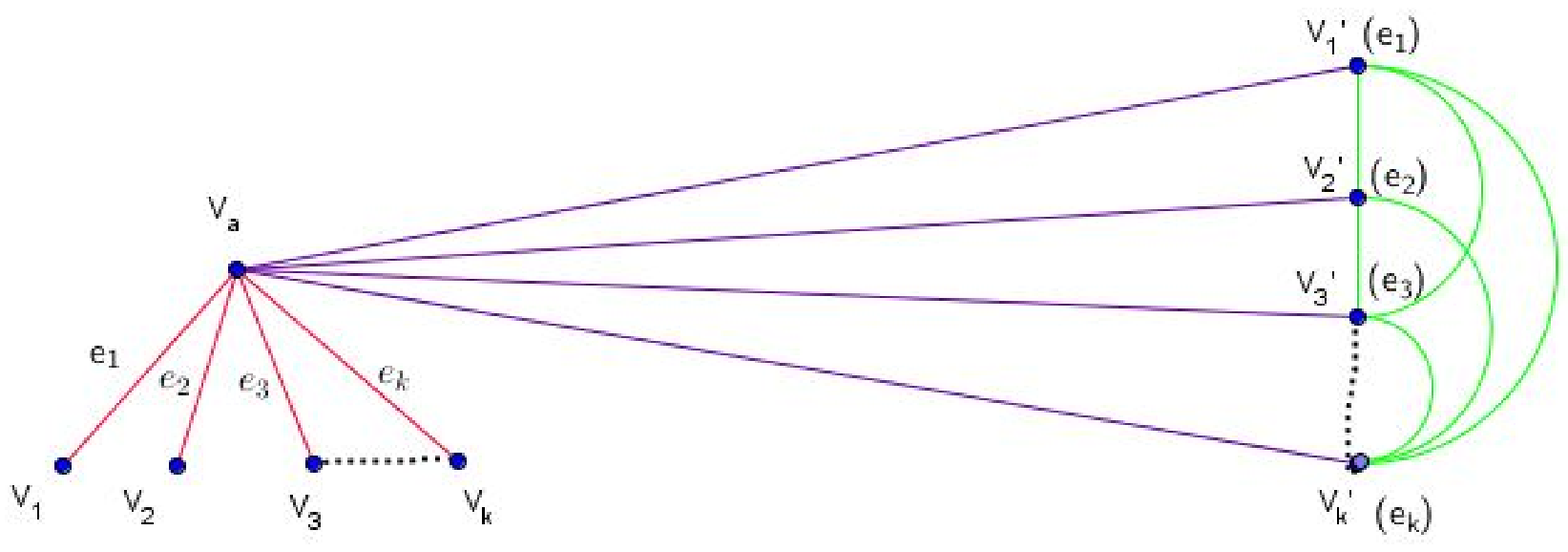} &
    \multicolumn{2}{c}{\includegraphics[height= 5cm, width=.50\textwidth]{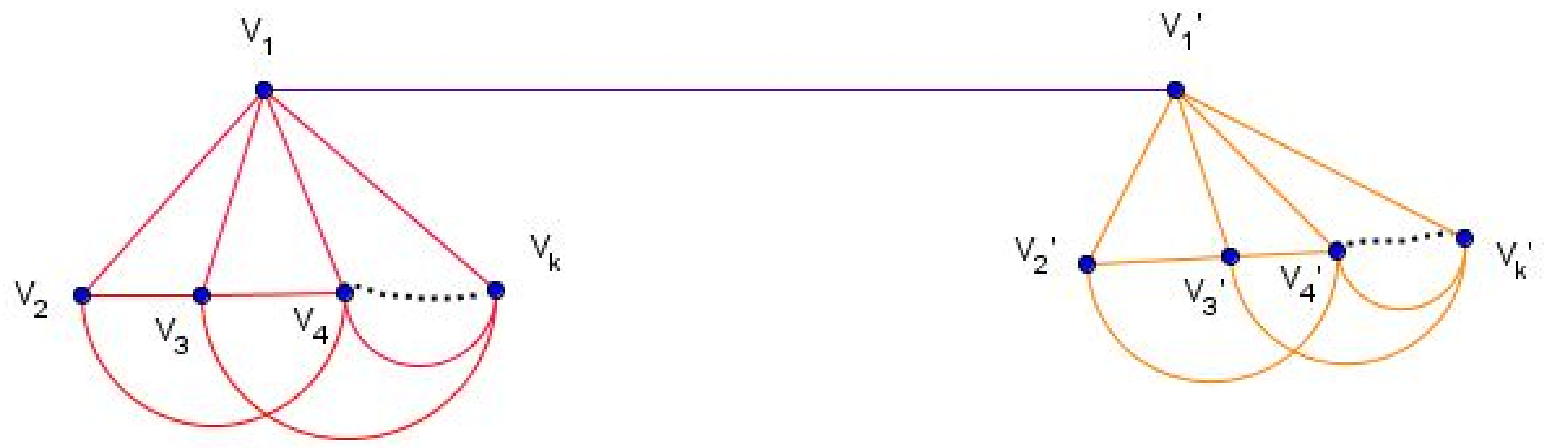}}
  \end{tabular}
  \caption{Characteristics of Vertex-Vertex and Edge-Vertex \label{fig:vertexedge}}
\end{figure}

\begin{theorems}\label{vertvertchar}
Given an arbitrary maximum degree vertex $v$, of degree $2k$, in a  total graph $H=T(G)$, ($G\ncong K_n$) it is a vertex vertex, if and only if its open neighbourhood satisfies the following properties:.
\begin{enumerate}
\item Its neighbours can be divided into two disjoint and exhaustive groups of $k$ vertices each, one corresponding to its vertex neighbours $N^{\mathcal{V}}(v)$ in the inverse total graph and the other corresponding to its incident edges $N^{\mathcal{E}}(v)$ in the inverse total graph.
\item The maximum degree of $G[N^{\mathcal{V}}(v) \Cup N^{\mathcal{E}}(v)]$ is $k$. This number is achieved by each vertex in $N^{\mathcal{E}}(v)$ and at least one vertex of $N^{\mathcal{V}}(v)$  falls short of this degree.
\end{enumerate}
\end{theorems}
\vspace{-0.3cm}
\begin{proof}
The statements follow from an inspection of Figure~\ref{fig:vertexedge}.
\vspace{-.7cm}\begin{flushright}$\Box$\end{flushright}
\end{proof}
\vspace{-0.3cm}

\begin{theorems}\label{edgevertchar}
Given an arbitrary maximum degree vertex $v$, of degree $2k$, in a  total graph $H=T(G)$, ($G\ncong K_n$) it is an edge vertex, if and only if  in the subgraph induced on its open neighbourhood:\\
Its neighbours can be partitioned into 2 maximal cliques of exactly $k$ vertices each consisting of one vertex vertex and $k-1$ edge vertices.
\end{theorems}
\vspace{-0.3cm}
\begin{proof}
The statement follows from an inspection of Figure~\ref{fig:vertexedge}.
\vspace{-.7cm}\begin{flushright}$\Box$\end{flushright}
\end{proof}
\vspace{-0.3cm}
In 4-regular candadate total graphs (total graphs of cycles), each vertex satisfies the conditions of both the above theorems. That is because the two parts of any valid partition are isomorphic to each other.
\begin{theorems}
Given any input (candidate) total graph $H=T(G)$, where $G\ncong K_n$ , it is indeed a total graph if and only if the graph $H'=T(G')$ is a total graph where $H'$ is obtained from $G$ by eliminating a vertex vertex of maximum degree along with all its edge neighbours in $H$.
\end{theorems}
\vspace{-0.3cm}
\begin{proof}
If the given graph is indeed a total graph, then it has a partition of its vertex set into vertex vertices and edge vertices. The algorithmic version of Theorem~\ref{vertvertchar} allows us to find a vertex vertex of maximum degree as well as identify its edge neighbours. Deleting the vertex along with its edge neighbours, effectively eliminates the vertex and the incident edges from the inverse total graph. We are left with the total graph of the graph with one vertex deleted. Recursing on the smaller graph, we obtain the partition or conclude that none exists if at some iteration a maximum degree vertex violates both Theorems~\ref{vertvertchar},~\ref{edgevertchar}. 
\vspace{-.7cm}\begin{flushright}$\Box$\end{flushright}
\end{proof}
\vspace{-0.1cm}
We thus have an algorithm which starts with a candidate total graph of a non-complete graph and decides whether it is indeed a total graph by recursing on the smaller graph or recourse to complete graph theorem.\\

{\bf \underline {Algorithm 2: Inverse Total Graph}}
\begin{enumerate}
\item Check if the given graph is the total graph of a complete graph using Algorithm~1. If so augument the vertices of that complete subgraph to the vertices obtained in earlier iterations and return. 
\item Else scan the vertices of maximum degree one by one, for conformity to Theorem~\ref{vertvertchar} until one, say $x$ is found. If during this process we find a vertex of maximum degree violating both Theorems~\ref{vertvertchar},~\ref{edgevertchar} then stop and conclude that it is not the total graph of any graph.
\item For the vertex $x$, from Step  2, partition its neighbours into vertex vertices and edge vertices by the algorithmic version of Theorem~\ref{vertvertchar}.
\item Add $x$ to the vertex set of the inverse total graph and repeat the steps with the graph obtained by deleting $x$ and its edge neighbours.
\item Return the set of vertices (and the induced subgraph on them) accumulated in Step 4 over all the iterations.
\end{enumerate}
\vspace{-0.3cm}
\section{Conclusions \& Future Directions}\label{SecConcl}
\vspace{-0.3cm}
We have proved properties of the vertex degrees of total graphs. We have developed a precise characterisation of the structure of the neighbourhoods of maximum degree vertices of the total graph of any graph. Combining these results we have designed an efficient iterative algorithm to compute the inverse total graph of a candidate total graph, or report that the graph is not a total graph. We also present a direct construction for the total graphs of paths, cycles and complete graphs. 

One interesting direction of future research is to see if a given $n$ and $m$ pair admits a connected unique total graph if any. One can also look at minimum number of dynamic graph operations (adding/deleting vertices and edges or moving edges around the graph) to transform a non-total graph into a total graph.
\vspace{-0.3cm}
\section{References}
\begin{enumerate}
	\item  M Behzad. A characterization of total graphs. Proceedings of the American Math-
ematical Society, 26(3):383-389, 1970.
\item  Mehdi Behzad. The total chromatic number of a graph: a survey. Combinatorial
	Mathematics and its Applications, pages 18, 1971.
\item  Mehdi Behzad and Heydar Radjavi. The total group of a graph. Proceedings of the
		American Mathematical Society, pages 158-163, 1968.
\item  Mehdi Behzad and Heydar Radjavi. Structure of regular total graphs. Journal of
			the London Mathematical Society, 1(1):433-436, 1969.
				\item  Robert L Hemminger and LowellWBeineke. Line graphs and line digraphs. Selected
				topics in graph theory, 1:291305, 1978.
					\item  J Krausz. Demonstration nouvelle dune theoreme de whitney sur les reseaux. Mat.
					Fiz. Lapok, 50(75-85):11, 1943.
					\item  ACMM van Rooij and H Wilf. The interchange graph of a nite graph. Acta
					Mathematica Hungarica, 16(3-4):263-269, 1965.
						\item  Vadim G Vizing. Some unsolved problems in graph theory. Russian Mathematical
	Surveys, 23(6):125141, 1968.
\end{enumerate} 

\end{document}